\documentclass[a4paper]{amsart}
\usepackage[margin=3.3cm]{geometry}
\usepackage[utf8]{inputenc}
\usepackage[english]{babel}
\usepackage{amssymb}
\usepackage{amsmath}
\usepackage{amsthm}
\usepackage{amsfonts}
\usepackage{array}
\usepackage{comment,stmaryrd}
\usepackage{dirtytalk}
\usepackage{enumitem}
\usepackage{hyperref}
\usepackage{mathrsfs}
\usepackage{tikz}
\usetikzlibrary{arrows,decorations.markings}
\usepackage{tikz-cd, mathtools}
\usepackage{mathabx}
\usepackage{xfrac}
\usepackage[aligntableaux=center]{ytableau}
\usepackage{todonotes}
\usepackage{graphicx}
\graphicspath{ {./Images/} }

\addto\captionsenglish{}

\theoremstyle{plain}
\newtheorem{theorem}{Theorem}[section]
\newtheorem{lemma}[theorem]{Lemma}

\newtheorem{proposition}[theorem]{Proposition}
\newtheorem{corollary}[theorem]{Corollary}

\theoremstyle{definition}
\newtheorem{definition}[theorem]{Definition}
\newtheorem{assumptions}[theorem]{Assumptions}
\newtheorem{example}[theorem]{Example}

\theoremstyle{remark}
\newtheorem{remark}[theorem]{Remark}

\newtheorem*{theorem*}{Theorem}
\newtheorem*{conjecture*}{Conjecture}

\DeclareMathOperator{\im}{Im}

\DeclareMathOperator{\End}{End}

\DeclareMathOperator{\Sym}{Sym}

\DeclareMathOperator{\Spec}{Spec}
\DeclareMathOperator{\Proj}{Proj}

\DeclareMathOperator{\Hilb}{Hilb}
\DeclareMathOperator{\Quot}{Quot}
\DeclareMathOperator{\NQV}{\mathfrak{M}}

\DeclareMathOperator{\Add}{Add}
\DeclareMathOperator{\MAdd}{\mathcal{A}\mathit{dd}}

\newcommand*{\SHom}{\mathcal{H}\kern -.5pt om}
\newcommand*{\SExt}{\mathcal{E}\kern -.5pt xt}

\DeclareMathOperator{\Aff}{\mathbb{A}}

\DeclareMathOperator{\GL}{GL}

\DeclareMathOperator{\SL}{SL}
\DeclareMathOperator{\Z}{\mathbb{Z}}

\DeclareMathOperator{\N}{\mathbb{N}}

\DeclareMathOperator{\C}{\mathbb{C}}
\DeclareMathOperator{\ida}{\mathfrak{a}}
\DeclareMathOperator{\Rep}{Rep}
\DeclareMathOperator{\Mat}{Mat}
\DeclareMathOperator{\Mo}{\mathcal{M}}
\DeclareMathOperator{\Sing}{\mathbb{A}^2\!/\Gamma}

\makeatletter
\let\@wraptoccontribs\wraptoccontribs
\makeatother

\title[Quiver Schemes for Nongeneric Stability and Cornering]{Quiver Schemes for Nongeneric Stability and Cornering}

\author{Lukas Bertsch}
\address{
Faculty of Mathematics, 
University of Vienna, 
Oskar-Morgenstern-Platz 1, 1090 Vienna, 
Austria}
\email{lukas.bertsch@univie.ac.at }

\author{Søren Gammelgaard}
\address{Dipartimento di Matematica "Federigo Enriques",
Università degli Studi di Milano,
Via Saldini 50 20123 Milano MI, Italy}
\email{soren.gammelgaard@unimi.it}


\begin{document}

\begin{abstract}
    We use the Le Bruyn--Procesi theorem to prove several results on quiver schemes for nongeneric stability conditions: We show that the stability-zero finite ADE-type Nakajima quiver schemes are reduced points and give an example of a closely related nonreduced quiver scheme. In broader generality, we prove that adding modules supported on stability-zero vertices induces closed embeddings of quiver schemes, and show how in many cases the quiver scheme associated with the cornered algebra defines a limit to this system of embeddings. As an application, we show that there is an isomorphism between the underlying reduced schemes of certain equivariant Quot schemes and Nakajima quiver varieties.
\end{abstract}

\maketitle

\section{Introduction}

In much recent work, Nakajima quiver schemes at nongeneric stability conditions were shown to be related to various moduli schemes associated to Kleinian (that is, canonical surface) singularities \cite{craw2024orbifoldquotschemesle, craw2019punctual, craw2026hilbertschemescrepantpartial, bertsch2026stackyresolutionskleiniansingularities}. A persistent issue when working with Nakajima quiver schemes in this context is that much of the established theory relies on point-wise arguments and does not address the possibility of non-reduced structure. In this paper, continuing on work of Craw-Yamagishi \cite{craw2026hilbertschemespointscanonical}, we use the Le Bruyn--Procesi theorem \cite[Theorem 1]{lebruynprocesi90} and its generalisation due to Lusztig \cite[Theorem 1.3]{Lusztig98} to establish multiple results on a large class of quiver GIT quotients, including Nakajima quiver schemes.

Nakajima quiver schemes or varieties, originally introduced by Nakajima \cite{nakajima1994instantons} and Kronheimer \cite{kronheimer_ale_89}, are moduli spaces of modules over a class of algebras called \emph{preprojective} algebras. Our first result is the following.

\begin{theorem}[Corollary \ref{cor:no-invariants-finite-type}]\label{thm:no-invariants-finite-type-intro}
        Let $\Pi$ be the (unframed) preprojective algebra of a possibly disconnected finite ADE Dynkin diagram. Let $K$ denote its vertex set and let $v = (v_i)_{i\in K} \in \N^K$. Then the Nakajima quiver scheme \[\Mo_{K,0}(\Pi, v) \coloneqq \Rep(\Pi, v) \sslash  \left(\prod_{i \in K} \GL(v_i)\right)\] is isomorphic to $\Spec(\C)$, i.e., a reduced point.
\end{theorem}

Here, the statement that $\Mo_{K,0}(\Pi, v)$ only has one point is well known, but that the scheme indeed is reduced, is new. We deduce it from an algebraic result of Malkin, Ostrik, and Vybornov \cite{MALKIN2006514}.

Next, we show that the representation scheme $\Rep(\Pi, v)$ itself is not reduced in some cases. Let $\delta_K \in \N^K$ be the restriction of the dimension vector corresponding to the minimal positive imaginary root on the affine Dynkin diagram of the same type as $\Pi$.

\begin{theorem}[Theorem \ref{thm:rep-not-reduced}]\label{thm:rep-not-reduced-intro}
    $\Rep(\Pi,\delta_K)$ is not reduced in types D and E.
\end{theorem} 

As far as we are aware, Theorem \ref{thm:rep-not-reduced-intro}  gives the first example of a nonreduced representation scheme associated to a preprojective algebra. On the other hand, it is still open whether this nonreducedness carries over to any Nakajima quiver variety that is obtained as a GIT quotient of $\Rep(\Pi,\delta_K)$ for some nontrivial stability condition.

For the proof of Theorem \ref{thm:rep-not-reduced-intro}, we identify the minimal resolution of a Kleinian singularity with the quiver variety $\Mo_{I,\zeta}(\Tilde{\Pi}, \delta)$, where $\Tilde{\Pi}$ is the \emph{affine} preprojective algebra, and then show (Theorem \ref{thm:kleinian-fibre}) that the scheme-theoretic fibre of this minimal resolution is isomorphic to $\Mo_{I,\zeta}(\Tilde{\Pi}^*, \delta)$, where $\Tilde{\Pi}^*$ is a quotient of $\Tilde{\Pi}$ by the arrows with head at the extended vertex. We also make use of a technique known as \emph{cornering}, which we will consider in broader generality in the second half of the paper.

We move on to study a larger class of algebras. Consider a quiver $Q$, and fix a decomposition of its vertex set \[Q_0 = F \sqcup J \sqcup K \; ,\] $J$ and $K$ non-empty, and write $H \coloneqq F \sqcup J$ and $I \coloneqq J \sqcup K$. Here, $F$ is the set of \emph{framing vertices}. We fix a stability vector $\zeta \in \Z^I$ such that
\begin{equation}\label{eq:nongen-stability-intro}
    \zeta|_K=0 \; .
\end{equation}
Let $\ida < \C Q$ be a homogenous ideal of relations and denote the quotient algebra by $A \coloneqq \C Q/\ida$. We assume that that there is an isomorphism \[A / (e_H) \cong \Pi\] of quiver algebras with vertex set $K$, where $e_H \in A$ is the idempotent of the set $H \subseteq Q_0$. The algebras $\Tilde{\Pi}$ and $\Tilde{\Pi}^*$ that appear in the proof of Theorem \ref{thm:kleinian-fibre} are examples of such $A$.

For any $v \in \N^{Q_0}$, we now consider the quiver scheme \[\Mo_{I,\zeta}(A,v) \coloneqq \Rep(A,v) \sslash_{\chi_{\zeta}} \prod_{i \in I} \GL(v_i)\] parametrising $\zeta$-semistable $A$-modules. Since stability vectors $\zeta$ of the form (\ref{eq:nongen-stability-intro}) are usually \emph{nongeneric}, meaning that there are $\zeta$-semistable $A$-modules that are not $\zeta$-stable, $\Mo_{I,\zeta}(A,v)$ will often not be a fine moduli space. A technique to get around this issue, used in e.g.\ \cite{craw2024orbifoldquotschemesle}, \cite{CGGS2}, \cite{craw2019punctual}, and \cite{bertsch2026stackyresolutionskleiniansingularities} is to replace the algebra $A$ with the \emph{cornered} algebra \[A_H \coloneqq e_H A e_H \; ,\] and the quiver scheme $\Mo_{I,\zeta}(A,v)$ with the scheme $\Mo_{J,\zeta|_J}(A_H,v|_H)$. In this paper, we develop this construction in greater generality. Under the assumptions stated above, we show that $A_H$ is finitely generated (Lemma \ref{lem:finite-generation-cornered-algebra}), and deduce that $A_H$ can also be written as the quotient of the path algebra of a quiver with vertex set $H$ (Corollary \ref{cor:finGenBetter}) so that we can once again construct moduli schemes $\Mo_{J,\zeta|_J}(A_H,v|_H)$. Then we prove the following theorem.

\begin{theorem}[Proposition \ref{prop:add-trivial}, Theorem \ref{thm:cornered-embedding}]\label{thm:cornered-embedding-intro}
    Let $v_H \in \N^H$. There is a system of closed embeddings of quiver schemes, indexed by the poset of vectors $v \in \N^{Q_0}$ satisfying $v|_H = v_H$: \[\cdots \hookrightarrow \Mo_{I,\zeta}(A,v) \hookrightarrow \Mo_{I,\zeta}(A,v') \hookrightarrow \cdots \hookrightarrow \Mo_{J,\zeta_J}(A_H,v_H)\] Here, $v,v' \in \N^{Q_0}$ are arbitrary vectors satisfying $v \leq v'$ and $v|_H = v'|_H = v_H$, and the closed embeddings on the left are induced by taking direct sums with vertex simple modules supported on $K$. The morphisms to $\Mo_{J,\zeta_J}(A_H,v_H)$ are induced by sending an $A$-module $V$ to the $A_H$-module $e_H V$.
    
    Furthermore, there exists $v \in \N^{Q_0}$ such that $v|_H = v_H$, and for any $v' \geq v$ with $v'|_H = v_H$, the embedding $\Mo_{I,\zeta}(A,v) \hookrightarrow \Mo_{J,\zeta_J}(A_H,v_H)$ is a bijection on closed points.
\end{theorem}

It is a compelling but generally open question whether the bijections on closed points appearing in Theorem \ref{thm:cornered-embedding-intro} are indeed isomorphisms of schemes (cf.\ Remark \ref{rem:is-bij-isomorphism}).

As an application, we prove in section \ref{sec:quot-nqv} the following theorem. It concerns the \emph{orbifold Quot schemes} $\Quot_{J}^{v_J} ([\Sing])$ introduced in \cite{CGGS2}, which are fine moduli spaces interpolating between the Hilbert schemes $\Hilb^n([\Sing])$ and $\Hilb^n(\Sing)$ associated to a finite subgroup $\Gamma<\SL(2,\C)$ and the corresponding Kleinian orbifold and singularity, respectively.

\begin{theorem}[Theorem \ref{main_thm_appendix}]\label{main_thm_appendix-intro}
    Let $J \subseteq I$ be a non-empty subset and $v_J \in \N^J$ be a dimension vector. The orbifold Quot scheme $\Quot_J^{v_J} ([\Sing])$ is non-empty if and only if the Nakajima quiver scheme  $\mathfrak{M}_{\zeta^{\bullet}}(v,\Lambda_0)$ is non-empty for some vector $v\in \N^I$ satisfying $v|_J=v_J$. In this case, for a sufficiently large choice of $v \in \N^I$ with $v|_J=v_J$, we have an isomorphism \[\Quot_J^{v_J} ([\Sing])_{\rm red} \cong  \mathfrak{M}_{\zeta^{\bullet}}(v,\Lambda_0)_{\rm red},\] where on both sides we take the reduced scheme structure.
\end{theorem}

Theorem \ref{main_thm_appendix-intro} was already stated in \cite{CGGS2}, but a gap in the proof of \cite[Lemma 4.3]{craw2019punctual}, on which \cite[Proposition 6.1]{CGGS2} builds, was pointed out to the authors by Yehao Zhou. Here, we bypass that gap. Another proof of Theorem \ref{main_thm_appendix-intro} appears in \cite{craw2024orbifoldquotschemesle}.

\subsection*{Acknowledgments}
We thank \'Ad\'am Gyenge and Bal\'azs Szendr\H{o}i for several conversations and comments. We thank Alastair Craw for his comments on a preliminary version of this paper, as well as for his explanations regarding \cite{craw2025lebruynprocesitheoremfollowing}. We also thank the anonymous reviewer for their helpful comments. While working on this, S. Gammelgaard was supported by the PRIN project 20223B5S8L `Birational geometry of moduli spaces and special varieties'.

\subsection*{Conventions}
We work throughout over $\C$. All quivers are assumed finite. Given a quiver $Q$, in its quiver algebra $\C Q$, we compose arrows in the same way as functions: for two paths $p,p'  \in \C Q$, the element $pp'$ denotes the path \emph{starting} with $p'$, ending with $p$ (where the composition makes sense).

For practical reasons we use the convention $0\in \N$. Given an index set $I$, we impose a partial order on the monoid $\N^I$, by saying that for $v, u\in \N^I$, we have $v\le u$ if $v_i\le u_i$ for every $i\in I$.

\section{Quiver Schemes and Le Bruyn--Procesi with Frozen Vertices}\label{sec:quiver-varieties}

Let $Q=(Q_0,Q_1,s,t)$ be a quiver, where $Q_0$ is a finite set of vertices, $Q_1$ a finite set of edges, and $s,t \colon Q_1 \to Q_0$ are the maps sending an edge to its \emph{source} and \emph{target} vertex. Let $\C Q$ be the path algebra of $Q$. $\C Q$ is naturally graded by \emph{path length}, which we denote by $l(p)$ for a homogeneous element $p \in \C Q$. The path algebra contains idempotents $e_i \in \C Q$, each corresponding to the paths of length zero at $i \in Q_0$, satisfying $1 = \sum_{i \in Q_0} e_i$. For any subset $I \subseteq Q_0$ we write $e_I = \sum_{i \in I} e_i$, which is also an idempotent.

Let $\ida$ be a two-sided ideal in $\C Q$, and denote by $A \coloneqq \C Q / \ida$ the quotient algebra. The images of the $e_i$ in $A$ are again idempotents which add up to $1$, and we will use the same notation for them. A finite-dimensional left $A$-module $V$ has a \emph{dimension vector} $\dim(V) \coloneqq (\dim V_i)_{i \in Q_0} \in \N^{Q_0}$, where $V_i \coloneqq e_i V$. The idempotents $e_i$ give $\C Q$ and $A$ the structure of $\C^{Q_0}$-algebras.  The pair $(Q, \ida)$ is also commonly known as a \emph{quiver with relations}. A representation of $(Q, \ida)$ is the same thing as a left module over $A$.

\begin{remark}\label{rem:terms}
    We use the following terminology: The word \emph{path} will denote an element of $\C Q$ given by a path in the underlying quiver. We use the term \emph{class} or \emph{class of a path} for an element of $A$ which is the image of a path under the map $\C Q\to A$. We use \emph{cycle} for a path of length at least 1 starting and ending at the same vertex, and \emph{cycle class} for the class of a cycle. Given a $\C Q$-module $V$, we say that $V$ is \emph{supported} on a vertex set $I \subset Q_0$ if $e_jV = 0$ for every vertex $j \in Q_0 \setminus I$. By a \emph{vertex simple} $\C Q$-module, we mean a one-dimensional module supported at a single vertex on which every arrow acts by zero.
\end{remark}

Given a dimension vector $v \in \N^{Q_0}$, there are affine representation schemes $\Rep(Q,v)$ and $\Rep(A,v)$ which act as parameter spaces for $\C Q$-, respectively $A$-module structures on a given finite-dimensional $\C^{Q_0}$-module $V = \bigoplus_{i \in Q_0} \C^{v_i}$. We have \[\Rep(Q,v) = \prod_{h \in Q_1} \Mat_{v_{t(h)} \times v_{s(h)}} \; ,\] where $\Mat_{k \times l}$ is the space of $(k \times l)$-matrices, considered as a scheme isomorphic to affine space $\Aff^{kl}$. Imposing the relations of $\ida$ gives a closed embedding \[\Rep(A,v) \hookrightarrow \Rep(Q,v) \; .\] 

The spaces $\Rep(Q,v)$ and $\Rep(A,v)$ do not quite parametrise isomorphism classes of modules, because the underlying vector space of the modules comes equipped with a basis. A change of basis is encoded by an element of the group \[\GL(v) \coloneqq \prod_{i \in Q_0} \GL(v_i) \; ,\] and we have compatible actions of $\GL(v)$ on $\Rep(Q,v)$ and $\Rep(A,v)$, respectively. In slightly broader generality, we may wish to change basis only at vertices in a subset $I \subset Q_0$, in which case we consider the restricted action of \[\GL(v|_I) = \prod_{i \in I} \GL(v_i) < \GL(v) \; .\] In this case, the vertices in $F \coloneqq Q_0 \setminus I$ are referred to as \emph{framing vertices}, and the action of $\GL(v|_I)$ encodes isomorphisms of \emph{framed modules}. \footnote{Craw and Yamagishi \cite{craw2025lebruynprocesitheoremfollowing} also call framing vertices \emph{frozen vertices}.}

\begin{remark}\label{rem:rep-functor}
    The representation scheme $\Rep(A,v)$ is a fine moduli space for families of algebra homomorphisms \[A \to \End \left( \bigoplus_{i \in Q_0} \C^{v_i} \right)\] that map $e_i$ to the projection onto the $i$-th direct summand. A group element $g \in \GL(v|_I)$ acts by composing such a homomorphism with the automorphism of $\End( \bigoplus_{i \in Q_0} \C^{v_i})$ given by conjugation with $g$.
\end{remark}

It is now natural to construct moduli spaces of (framed) modules as GIT quotients by this action. For this, pick a stability vector $\zeta \in \Z^I$, and let \[\chi_{\zeta} \colon \GL(v|_I) \rightarrow \C^{\times} \; , \quad (g_i)_{i \in I} \mapsto \prod_{i \in I} \det(g_i)^{\zeta_i} \; \] be the corresponding character. Define the \emph{quiver schemes} as GIT quotients \footnote{These schemes are often called quiver \emph{varieties}, but this can be misleading as they are not necessarily reduced or irreducible. Occasionally, this is handled by always taking the underlying reduced scheme of these GIT quotients, but we use the word `scheme' here to indicate that we are actually working with the GIT quotients.} \[\Mo_{I,\zeta}(Q,v) \coloneqq \Rep(Q,v) \sslash_{\chi_{\zeta}} \GL(v|_I) \; , \quad \text{and} \quad \Mo_{I,\zeta}(A,v) \coloneqq \Rep(A,v) \sslash_{\chi_{\zeta}} \GL(v|_I) \; .\] Note that, since $\GL(v|_I)$ is reductive, the closed embedding of representation schemes carries over to a closed embedding of quotients \[\Mo_{I,\zeta}(A,v) \hookrightarrow \Mo_{I,\zeta}(Q,v) \; .\]

Recall furthermore that by construction, $\Mo_{I,\zeta}(Q,v)$ is projective over the affine scheme \[\Mo_{I,0}(Q,v) = \Spec \left(\C[\Rep(Q,v)]^{\GL(v|_I)} \right) \; ,\] and analogously with $A$ in the place of $Q$. The following theorem, which is a generalisation of the Le Bruyn--Procesi Theorem \cite[Theorem 1]{lebruynprocesi90} due to Lusztig \cite[Theorem 1.3]{Lusztig98} describes this affine quotient (see also \cite[Theorem 1.2]{craw2025lebruynprocesitheoremfollowing}).

\begin{theorem}[Le Bruyn--Procesi with frozen vertices]\label{Thm:LBP-frozen}
    The invariant ring \[\C[\Rep(Q,v)]^{\GL(v|_I)}\] is generated by
    \begin{enumerate}
        \item \label{Thm:LBP-frozen:traces} traces of products of matrices along cycles in $Q$ which only pass through vertices in $I$, and 
        \item \label{Thm:LBP-frozen:entries} entries of products of matrices along paths which start and end outside of $I$.
    \end{enumerate}
    Furthermore, the equivariant surjection of coordinate rings $\C[\Rep(Q,v)] \rightarrow \C[\Rep(A,v)]$ induces a surjection \[\C[\Rep(Q,v)]^{\GL(v|_I)} \rightarrow \C[\Rep(A,v)]^{\GL(v|_I)} \; ,\] so $\C[\Rep(A,v)]^{\GL(v|_I)}$ as well is generated by images of elements of the form (\ref{Thm:LBP-frozen:traces}) and (\ref{Thm:LBP-frozen:entries}).
\end{theorem}

\section{The finite preprojective Algebra}\label{sec:fin-preproj}

\subsection{Reduced and nonreduced quiver schemes}

In this section we consider the preprojective algebra $\Pi$ of a finite, possibly disconnected ADE Dynkin diagram with vertex set $K$. Recall that $\Pi$ is the quotient of the path algebra of the doubled quiver associated with the Dynkin diagram by the homogeneous ideal generated by the preprojective relations (for a precise definition, see e.g.\ \cite[Definition 5.2]{kirillov2016quiver}). 

The following result has been known since the discovery of preprojective algebras. \cite{gelfand1979model}

\begin{theorem}\label{thm:finite-dimensional-finite-type}
    $\Pi$ is finite-dimensional.
\end{theorem}

In particular, there is an upper bound to the length of an element of $\Pi$. Similarly, we have the following result, shown in \cite[Theorem 1.4a]{MALKIN2006514}, which implies that cycle classes in $\Pi$ of length $\geq 1$, up to cyclic permutation, are always zero.

\begin{theorem}[Malkin--Ostrik--Vybornov]\label{thm:cocenter-finite-type}
    \[\frac{\Pi}{[\Pi,\Pi]} \cong \C^K \; ,\] where the idempotents on the right-hand side correspond to the images of the vertex idempotents in $\Pi$.
\end{theorem}

Moduli spaces of $\Pi$-modules, as constructed in Section \ref{sec:quiver-varieties}, have been studied extensively. Here we consider the affine schemes $\Mo_{K,0}(\Pi, v)$ parametrising unframed modules with the stability condition $\zeta =0$. It is well known that the algebraic set underlying this scheme consists of only one point \cite[Prop.\ 1.2, Ex.\ 1.3]{crawley2002decomposition}, \cite[Thm.\ 10.16]{kirillov2016quiver}. As a consequence of Theorem \ref{thm:cocenter-finite-type}, we get a refinement of this statement, providing also reducedness of the quotient as a scheme.

\begin{corollary}\label{cor:no-invariants-finite-type}
    Let $\Pi$ be the (unframed) preprojective algebra of a possibly disconnected finite Dynkin diagram. Let $K$ denote its vertex set and let $v \in \N^K$. Then \[\Mo_{K,0}(\Pi, v) = \Rep(\Pi, v) \sslash \GL(v)\] is isomorphic to $\Spec(\C)$.
\end{corollary}
\begin{proof}
    The Le Bruyn--Procesi Theorem states that $\C[\Rep(\Pi, v)]^{\GL(v)}$ is generated by traces of cycle classes in $\Pi$. Such traces are functions of the images of these cycle classes in $\Pi/[\Pi,\Pi]$, so $\C[\Rep(\Pi, v)]^{\GL(v)}$ is generated by the image of the linear map \[\frac{\Pi}{[\Pi,\Pi]} \to \C[\Rep(\Pi, v)]^{\GL(v)}\] which sends the class of a cycle in $\Pi$ to the trace function of that cycle. But $\Pi/[\Pi,\Pi] = \prod_{i \in K} \C e_i$ by Theorem \ref{thm:cocenter-finite-type} and all idempotents $e_i$ for $i \in K$ are mapped to constants under this map, so we have \[\C[\Rep(\Pi, v)]^{\GL(v)} = \C \; .\]
\end{proof}

Consider now a connected, finite ADE Dynkin diagram with vertex set $K$, along with the affine Dynkin diagram of the same type having vertex set $\{0\} \sqcup K$. Let $\delta \in \N^{\{0\} \sqcup K}$ be the dimension vector corresponding to the minimal positive imaginary root on the affine Dynkin diagram, and let $\delta_K = \delta|_K$ be its restriction to the finite Dynkin diagram. Let $\Pi$ and $\tilde{\Pi}$ be the associated finite and affine preprojective algebras, respectively. We have the following result, which we prove in section \ref{sec:proof-rep-not-reduced}.

\begin{theorem}\label{thm:rep-not-reduced} 
    $\Rep(\Pi,\delta_K)$ is not reduced when the Dynkin diagram is of type $D_n (n \geq 4)$, $E_6$, $E_7$, or $E_8$.
\end{theorem}

The representation schemes $\Rep(\Pi,\delta_K)$ are special instances of the representation schemes $\Rep(\tilde{\Pi},v)$ of the affine preprojective algebra $\tilde{\Pi}$, obtained by choosing $v=(0,\delta_K)$. It follows from a Theorem of Gan--Ginzburg \cite[Theorem 1.6]{ganginzburg2006} that the latter are reduced for $v=n \delta$, $n \in \N$. To our knowledge, whether reducedness holds for arbitrary dimension vectors has been an open problem, to which Theorem \ref{thm:rep-not-reduced} now provides a counterexample. However, it remains unknown whether the Nakajima quiver schemes given as GIT quotients of the representation spaces $\Rep(\tilde{\Pi}, v)$ for nonzero stability can be nonreduced.

\begin{example}
    Before we prove Theorem \ref{thm:rep-not-reduced} in full generality, we give a hands-on proof of nonreducedness in type $D_4$. In this case, the representation scheme $\Rep(\Pi,\delta_K)$ can be presented via $(1\times2)$-matrix-valued coordinates $a,b,c$ and $(2\times1)$-matrix-valued coordinates $d,e,f$, subject to the preprojective relations (up to a choice of signs) \[da=eb+fc \; , \; ad=0 \; , \; be=0 \; , \; \text{and} \; cf = 0 \; . \]

    \begin{center}\begin{tikzpicture}[ivert/.style={circle,draw=black,thick,inner sep=0pt,minimum size=5mm}]
    \node (bl) at (-1.2,-.75) [ivert] {};
    \node (br) at (1.2,-.75) [ivert] {};
    \node (c) at (0,0) [ivert] {};
    \node (t) at (0,1.4) [ivert] {};

    \draw[->] (bl) to ["$e$", bend left=20] (c) ;
    \draw[->] (c) to ["$b$", bend left=20] (bl);
    \draw[->] (br) to ["$f$", bend left=20] (c);
    \draw[->] (c) to ["$c$", bend left=20] (br);
    \draw[->] (t) to ["$d$", bend left=20] (c);
    \draw[->] (c) to ["$a$", bend left=20] (t);
    \end{tikzpicture}\end{center}

    To show that $\Rep(\Pi,\delta_K)$ is not reduced, it suffices to show that there is a nonzero regular function on $\Rep(\Pi,\delta_K)$ that vanishes at every $\C$-valued point. We will show that \[ebfc \colon \Rep(\Pi,\delta_K) \longrightarrow \Mat_{2 \times 2}\] is such a function.

    First, we show that $ebfc$ vanishes at every $\C$-valued point. Let any such point be defined by matrices $A,B,C \in \Mat_{1 \times 2}(\C)$, $D,E,F \in \Mat_{2 \times 1}(\C)$, defining the values of the corresponding lower-case coordinates at that point. We need to show that $EBFC = 0$. This is clear if either $EB=0$ or $FC=0$, so we will assume these to be nonzero. Furthermore, if $DA=0$, we would have $EB=-FC$, and so $EBFC = -EBEB = 0$ since $BE=0$. Hence, we also assume that $DA$ is nonzero.

    By these assumptions, $DA$, $EB$ and $FC$ are all $(2 \times 2)$-matrices of rank $1$. Furthermore, they are all square-zero, which implies that their kernels and images coincide. We will show that $EB$ and $FC$ have the same image. Indeed, suppose that $x \in \im(EB)$ and $y \in \im(FC)$ are linearly independent. Then \[DAx = EBx+FCx = FCx \quad \text{and} \quad DAy = EBy+FCy = EBy \; .\] Since $x$ spans the image and kernel of $EB$ and $y$ is linearly independent of $x$, $DAy = EBy$ is a nonzero scalar multiple of $x$. The same argument shows that $DAx = FCx$ is a nonzero scalar multiple of $y$. But since $x$ and $y$ were assumed to be linearly independent, this would imply that $DA$ has rank $2$, a contradiction. We conclude that $EB$ and $FC$ have the same image, and hence, the same kernel. This implies that they are nonzero scalar multiples of one another, say $EB = \lambda FC$, and hence, \[EBFC = \lambda EBEB = 0 \; .\]

    It remains to show that $ebfc$ is not the zero function. For this, consider the $\C[\epsilon]/(\epsilon^2)$-valued point of $\Rep(\Pi,\delta_K)$ defined by 
    \[A= \begin{pmatrix}
        0 & 2
    \end{pmatrix} \; , \;
    B= \begin{pmatrix}
        \epsilon & 1
    \end{pmatrix} \; , \;
    C= \begin{pmatrix}
        -\epsilon & 1
    \end{pmatrix} \; , \;
    D = \begin{pmatrix}
        1 \\ 0
    \end{pmatrix} \; , \;
    E = \begin{pmatrix}
        1 \\ -\epsilon
    \end{pmatrix} \; , \;
    F = \begin{pmatrix}
        1 \\ \epsilon
    \end{pmatrix} \; .\]
    One easily verifies that these matrices satisfy the preprojective relations, but that 
    \[EBFC = \begin{pmatrix}
        0 & 2\epsilon \\ 0 & 0
    \end{pmatrix} \neq 0 \; .\]
    
\end{example}

\subsection{The minimal resolution and its central fibre}\label{sec:proof-rep-not-reduced}

We denote by $\Tilde{\Pi}$ the preprojective algebra associated to the affine ADE Dynkin diagram of the same type as $\Pi$. We let $I$ be the vertex set of the affine diagram, and denote the added vertex by $0 \in I$, so that $K = I \setminus \{0\}$. Then \[\frac{\Tilde{\Pi}}{(e_0)} = \Pi \; . \] Let $\Tilde{\Pi}_0 \coloneqq e_0 \Tilde{\Pi} e_0$ be the algebra consisting of the classes of paths in $\Tilde{\Pi}$ starting and ending at $0$. (This is called the \emph{cornered} algebra, a notion that we will study in greater detail in section \ref{sec:cornered}.) There is a well-known isomorphism \[\Tilde{\Pi}_0 \cong \C[x,y]^{\Gamma} \; ,\] where the right-hand side is the coordinate ring of a Kleinian singularity for a finite subgroup $\Gamma < \SL(2, \C)$ \cite{reiten_van-den-bergh_two-dimensional-tame_89} (cf.\ \cite[Theorem 0.1]{crawley1998noncommutative}). Viewing $\Tilde{\Pi}_0$ as a quiver algebra with vertex set $\{0\}$, the quiver GIT quotient constructed in section \ref{sec:quiver-varieties} is precisely this Kleinian singularity: \[\Mo_{\{0\},0}(\Tilde{\Pi}_0, 1) = \Rep(\Tilde{\Pi}_0, 1) \cong \Spec(\Tilde{\Pi}_0) \cong \Sing \; .\] (Note that the gauge group $\GL(1)$ acts trivially in this case.)

Let $\delta \in \N^I$ be the dimension vector corresponding to the minimal positive imaginary root. Note that $\delta_0 = 1$. Furthermore, let $\zeta \in \Z^I$ be defined as follows. \[\zeta_i \coloneqq \begin{cases} 1 \; , & \text{when} \; i \neq 0 \; ; \\ -\sum_{i \in K} \delta_i \; , & \text{when} \; i = 0 \; .\end{cases}\] We construct the GIT quotient $\Mo_{I,\zeta}(\Tilde{\Pi}, \delta)$ as in section \ref{sec:quiver-varieties} (without framing vertices).

Furthermore, we consider the algebra $\Tilde{\Pi}^*$, which is the quotient of $\Tilde{\Pi}$ by the images of all arrows going from $K$ to $0$. The following illustrations show the quivers generating the algebras $\Tilde{\Pi}$, $\Tilde{\Pi}^*$ and $\Pi$ in type $D_5$.

\vspace{6mm}

\begin{center}
\begin{tikzpicture}[framing/.style={rectangle,draw=black,thick,inner sep=0pt,minimum size=5mm},ivert/.style={circle,draw=black,thick,inner sep=0pt,minimum size=5mm}]
    \node (a) at (0,0) {$\Tilde{\Pi}$};
    \node (tl) at (2,.8) [ivert] {$0$};
    \node (bl) at (2,-.8) [ivert] {};
    \node (l) at (3,0) [ivert] {};
    \node (r) at (4.3,0) [ivert] {};
    \node (tr) at (5.3,.8) [ivert] {};
    \node (br) at (5.3,-.8) [ivert] {};

    \draw[->] (tl) to [bend left=20] (l);
    \draw[->] (l) to [bend left=20] (tl);
    \draw[->] (bl) to [bend left=20] (l);
    \draw[->] (l) to [bend left=20] (bl);
    \draw[->] (l) to [bend left=20] (r);
    \draw[->] (r) to [bend left=20] (l);
    \draw[->] (tr) to [bend left=20] (r);
    \draw[->] (r) to [bend left=20] (tr);
    \draw[->] (br) to [bend left=20] (r);
    \draw[->] (r) to [bend left=20] (br);
\end{tikzpicture}
\\ \vspace{6mm}
\begin{tikzpicture}[framing/.style={rectangle,draw=black,thick,inner sep=0pt,minimum size=5mm},ivert/.style={circle,draw=black,thick,inner sep=0pt,minimum size=5mm}]
    \node (a) at (0,0) {$\Tilde{\Pi}^*$};
    \node (tl) at (2,.8) [ivert] {$0$};
    \node (bl) at (2,-.8) [ivert] {};
    \node (l) at (3,0) [ivert] {};
    \node (r) at (4.3,0) [ivert] {};
    \node (tr) at (5.3,.8) [ivert] {};
    \node (br) at (5.3,-.8) [ivert] {};

    \draw[->] (tl) to (l);
    \draw[->] (bl) to [bend left=20] (l);
    \draw[->] (l) to [bend left=20] (bl);
    \draw[->] (l) to [bend left=20] (r);
    \draw[->] (r) to [bend left=20] (l);
    \draw[->] (tr) to [bend left=20] (r);
    \draw[->] (r) to [bend left=20] (tr);
    \draw[->] (br) to [bend left=20] (r);
    \draw[->] (r) to [bend left=20] (br);
\end{tikzpicture}
\\ \vspace{6mm}
\begin{tikzpicture}[framing/.style={rectangle,draw=black,thick,inner sep=0pt,minimum size=5mm},ivert/.style={circle,draw=black,thick,inner sep=0pt,minimum size=5mm}]
    \node (a) at (0,0) {$\Pi$};
    \node (bl) at (2,-.8) [ivert] {};
    \node (l) at (3,0) [ivert] {};
    \node (r) at (4.3,0) [ivert] {};
    \node (tr) at (5.3,.8) [ivert] {};
    \node (br) at (5.3,-.8) [ivert] {};

    \draw[->] (bl) to [bend left=20] (l);
    \draw[->] (l) to [bend left=20] (bl);
    \draw[->] (l) to [bend left=20] (r);
    \draw[->] (r) to [bend left=20] (l);
    \draw[->] (tr) to [bend left=20] (r);
    \draw[->] (r) to [bend left=20] (tr);
    \draw[->] (br) to [bend left=20] (r);
    \draw[->] (r) to [bend left=20] (br);
\end{tikzpicture}
\end{center}

\vspace{6mm}

\begin{theorem}\label{thm:kleinian-fibre}
    We have the following Cartesian diagram of schemes.
    \[\begin{tikzcd}
        \Mo_{I,\zeta}(\Tilde{\Pi}^*, \delta) \arrow[r] \arrow[d] & \Mo_{I,\zeta}(\Tilde{\Pi}, \delta) \arrow[d, "\pi"] \\
        * \arrow[r] & \Mo_{\{0\},0}(\Tilde{\Pi}_0, 1)
    \end{tikzcd}\]
    Here, the morphism $\pi$, induced by sending a $\Tilde{\Pi}$-module $V$ to its $0$-component $V_0 = e_0 V$, is the minimal resolution of the Kleinian singularity, the morphism at the bottom is the embedding of the singular point into $\Mo_{\{0\},0}(\Tilde{\Pi}_0)$, and the morphism on top is the closed embedding induced by the quotient homomorphism $\Tilde{\Pi} \to \Tilde{\Pi}^*$.
\end{theorem}
\begin{proof}
    First, we construct the morphism $\pi$. Sending a $\Tilde{\Pi}$-module $V$ to its $0$-component $V_0$ induces a map of representation schemes $\Rep(\Tilde{\Pi}, \delta) \to \Rep(\Tilde{\Pi}_0,1)$ (cf.\ proof of Theorem \ref{thm:cornered-embedding}(\ref{thm:cornered-embedding:1}) below). This morphism is invariant under the action of $\GL(\delta)$ on $\Rep(\Tilde{\Pi}, \delta)$. This defines $\pi$ and shows that $\pi$ factors as $\pi = \pi_0 \pi'$ through the affine quotient \[\Mo_{I,\zeta}(\Tilde{\Pi}, \delta) \xlongrightarrow{\pi'} \Mo_{I,0}(\Tilde{\Pi}, \delta) \xlongrightarrow{\pi_0} \Mo_{\{0\},0}(\Tilde{\Pi}_0, 1) \; .\]
    
    By a result of Kronheimer (\cite{kronheimer_ale_89}, see also \cite[Theorem 2.3]{craw2026hilbertschemescrepantpartial}), the VGIT morphism $\pi'$ is the minimal resolution of the same Kleinian singularity; in particular, there are isomorphisms \[\Mo_{I,0}(\Tilde{\Pi}, \delta) \cong \Sing \cong \Mo_{\{0\},0}(\Tilde{\Pi}_0, 1) \; .\] We claim that $\pi_0$ is also an isomorphism. Indeed, both are affine varieties whose coordinate rings are generated by traces of cycles of the form (\ref{Thm:LBP-frozen:traces}) in Theorem \ref{Thm:LBP-frozen}. More precisely, the coordinate ring of $\Mo_{I,0}(\Tilde{\Pi}, \delta)$ is generated by traces of cycles in $\Tilde{\Pi}$, and the coordinate ring of $\Mo_{\{0\},0}(\Tilde{\Pi}_0, 1)$ is generated by traces of cycles in $\Tilde{\Pi}$ passing through $0$. It is a consequence of Theorem \ref{thm:cocenter-finite-type} that the trace of any cycle of length $\geq 1$ is equal to a linear combination of traces of cycles passing through $0$. This implies that the coordinate ring of $\Mo_{\{0\},0}(\Tilde{\Pi}_0, 1)$ surjects onto the coordinate ring of $\Mo_{I,0}(\Tilde{\Pi}, \delta)$ along $\pi_0$ (this argument is given in broader generality and more detail in the proof of Theorem \ref{thm:cornered-embedding}(\ref{thm:cornered-embedding:1})). Since $\Mo_{I,0}(\Tilde{\Pi}, \delta)$ and $\Mo_{\{0\},0}(\Tilde{\Pi}_0, 1)$ were isomorphic to begin with, this surjection is an isomorphism by the well-known fact that a surjective homomorphism from a finitely generated commutative algebra to itself is an isomorphism.

    By King's analysis of GIT stability for quivers \cite[Proposition 5.3]{king1994moduli}, $\Mo_{I,\zeta}(\Tilde{\Pi}, \delta)$ is a fine moduli space for zero-generated $\Tilde{\Pi}$-modules of dimension vector $\delta$. Analogously, it is clear that $\Mo_{\{0\},0}(\Tilde{\Pi}_0, 1)$ is a fine moduli space for $\Tilde{\Pi}_0$-modules of dimension $1$. The singular point $*$ in $\Mo_{\{0\},0}(\Tilde{\Pi}_0, 1)$ represents the unique $\Tilde{\Pi}_0$-module on which any path of length $\geq 1$ acts as zero. Hence, its scheme-theoretic fibre $\pi^{-1}(*) \subset \Mo_{I,0}(\Tilde{\Pi}, \delta)$ is a fine moduli space for those $0$-generated $\Tilde{\Pi}$-modules $V$ such that all paths of length $\geq 1$ starting and ending at $0$ act as zero on $V_0$. $\Mo_{I,\zeta}(\Tilde{\Pi}^*, \delta)$, on the other hand, is a fine moduli space for those $0$-generated $\Tilde{\Pi}$-modules on which all edges going from $K$ to $0$ act by zero.

    We now show that the subfunctors represented by $\Mo_{I,\zeta}(\Tilde{\Pi}^*, \delta)$ and $\pi^{-1}(*)$ in $\Mo_{I,\zeta}(\Tilde{\Pi}, \delta)$ are the same. Let $V$ be a flat family of $0$-generated $\Tilde{\Pi}$-modules over a test scheme $T$. The condition of being $0$-generated is equivalent to the condition that the natural homomorphism of quasi-coherent sheaves of $\Tilde{\Pi}$-modules on $T$ \[e_i \Tilde{\Pi} e_0 \otimes V_0 \longrightarrow V_i\] is a surjection for every $i \in I$.
    
    Assume that $V$ represents a morphism $T \to \Mo_{I,\zeta}(\Tilde{\Pi}^*, \delta)$, i.e., any edge going from $K$ to $0$ acts as the zero endomorphism on $V$. Since the class of any path of length $\geq 1$ in $\Tilde{\Pi}$ starting and ending in $0$ has to pass through $K$, any such path will act by zero on $V$, showing that $V$ also represents a morphism $T \to \pi^{-1}(*)$.
    
    Conversely, assume that $V$ represents a morphism $T \to \pi^{-1}(*)$. Let $h$ be an edge going from a vertex $i \in K$ to $0$. We have the following homomorphisms of quasi-coherent sheaves, where we use $h$ to indicate the morphism induced by the edge $h$. \[e_i \Tilde{\Pi} e_0 \otimes V_0 \longrightarrow V_i \xlongrightarrow{h} V_0 \; .\] Since $V$ is $0$-generated, the left-hand homomorphism here is surjective, but since any path starting and ending in $0$ that also passes through $i$ acts by $0$, the composition of the two homomorphisms is zero. This implies that $h$ must be the zero morphism. Hence, $V$ also represents a morphism $T \to \Mo_{I,\zeta}(\Tilde{\Pi}^*, \delta)$.
\end{proof}

Finally, we use Theorem \ref{thm:kleinian-fibre} to prove Theorem \ref{thm:rep-not-reduced}.

\begin{proof}[Proof of Theorem \ref{thm:rep-not-reduced}]
    According to Theorem \ref{thm:kleinian-fibre}, $\Mo_{I,\zeta}(\Tilde{\Pi}^*, \delta)$ is isomorphic to the scheme-theoretic fibre of the minimal resolution of a Kleinian singularity.
    
    This fibre is non-reduced in types $D$ and $E$. Indeed, it was shown by Artin \cite{Artin66IsolatedRational} that the fundamental cycle of the fibre is equal to the unique minimal positive cycle supported on the fibre intersecting each of its irreducible components non-positively. The coefficients of the components of the fibre in this cycle are given by the vector $\delta_K$, which has entries $> 1$ in any type $D_n$, $E_6$, $E_7$, $E_8$. Since any GIT quotient of a reduced scheme is reduced, this implies that $\Rep(\Tilde{\Pi}^*,\delta)$ is not reduced.

    Note that $\Tilde{\Pi}^*$ is the finite preprojective algebra with an additionally adjoined vertex $0$ and one (in types $D$ and $E$) or two (in type $A$) arrows from $0$ to $K$. Since these additional arrows are not affected by any relations, and since $\sum_{0 \to i} \delta_i \delta_0 = 2$, we have \[\Rep(\Tilde{\Pi}^*,\delta) \cong \Rep(\Pi,\delta_K) \times \Aff^2\] in all types. In particular, $\Rep(\Pi,\delta_K)$ is also not reduced in types $D$ and $E$.
\end{proof}

\section{Quiver Schemes at Nongeneric Stability Conditions}\label{sec:special-stability}

In what follows, we describe applications of Theorems \ref{Thm:LBP-frozen}, \ref{thm:finite-dimensional-finite-type}, and \ref{thm:cocenter-finite-type} to the study of more general quiver schemes at special stability conditions. For now, we fix the following assumptions and notation.

\begin{assumptions}\label{asn:1}
    Let $Q=(Q_0,Q_1,s,t)$ be a quiver, $\mathfrak{a} < \C Q$ a two-sided ideal, and $A = \C Q/\mathfrak{a}$ the quotient algebra. Let $Q_0 = F \sqcup J \sqcup K \; ,$ where $K$ is non-empty, and let $I \coloneqq J \sqcup K$ and $H \coloneqq F \sqcup J$. Fix a vector $\zeta \in \Z^I$ so that \[\zeta_i = 0 \quad \text{for} \quad i \in K \; ,\] and write $\zeta_J \coloneqq \zeta|_J$.
\end{assumptions}

The following lemma is foundational to our results (cf.\ \cite[Proof of Theorem 4.7]{craw2026hilbertschemespointscanonical}).

\begin{lemma}\label{lem:successive-GIT}
    Under Assumptions \ref{asn:1}, the quiver schemes $\Mo_{I,\zeta}(A,v)$ can equivalently be constructed as
    \[\Mo_{I,\zeta}(A,v) = \left(\Rep(A,v) \sslash \GL(v|_K) \right) \sslash_{\chi_{\zeta_J}} \GL(v|_J) \; ,\]
    where first we take the affine GIT quotient \[\Rep(A,v) \sslash \GL(v|_K) = \Spec \left(\C[\Rep(A,v)]^{\GL(v|_K)} \right) \; ,\] and then the projective GIT quotient with respect to the restricted stability vector $\zeta_J$.
\end{lemma}
\begin{proof}
    We write $\chi = \chi_{\zeta}$ and $\chi_J = \chi_{\zeta_J} = \chi|_{\GL(v|_J)}$. By definition, the quiver scheme is given as the projective GIT quotient \[\Mo_{I,\zeta}(A,v) = \Proj \left(\bigoplus_{k \geq 0}\C[\Rep(A,v)]^{\GL(v|_I),\chi^k} \right) \; ,\] where \[\C[\Rep(A,v)]^{\GL(v|_I),\chi^k} = \{ f \in \C[\Rep(A,v)] \; | \; g f = \chi(g)^k f \} \; .\] Now, by the assumption on $\zeta$, we have $\chi(g) = 1$ for $g \in \GL(v|_K)$ and $\chi(g) = \chi_J(g)$ for $g \in \GL(v|_J)$. Therefore, the condition $g f = \chi(g)^k f$ for all $g \in \GL(v|_I)$ is equivalent to the two conditions
    \begin{itemize}
        \item $g f = f$ for $g \in \GL(v|_K)$, i.e.\ $f \in \C[\Rep(A,v)]^{\GL(v|_K)}$, and
        \item $g f = \chi_J(g)^k f$ for $g \in \GL(v|_J)$.
    \end{itemize}
    In other words, we have \[\C[\Rep(A,v)]^{\GL(v|_I),\chi^k} = \left(\C[\Rep(A,v)]^{\GL(v|_K)}\right)^{\GL(v|_J),(\chi_J)^k} \; ,\] and the statement follows.
\end{proof}

We now show that taking a direct sum with a module $V_K$ that is supported only on $K$ gives a well-defined embedding of quiver schemes for stability $\zeta$. Later, we will apply this result to the case where $V_K$ is a  direct sum of vertex simple modules.

\begin{proposition}\label{prop:add-trivial}
     Consider a module $V_K$ of $A$ that is supported on $K$, i.e., $e_i V_K = 0$ for $i \in H$. Let $v \in \N^{Q_0}$ and $\hat{v} \coloneqq v + \dim(V_K)$. Then the operation \[V \mapsto V \oplus V_K \; ,\] where $V$ is any $v$-dimensional $A$-module, induces $\GL(v|_J)$-equivariant closed embeddings of affine quotients \[\Add_{V_K} \colon \Rep(A,v) \sslash \GL(v|_K) \hookrightarrow \Rep(A,\hat{v}) \sslash \GL(\hat{v}|_K)\] and closed embeddings of quiver schemes \[\MAdd_{V_K} \colon \Mo_{I,\zeta}(A,v) \hookrightarrow \Mo_{I,\zeta}(A,\hat{v}) \; .\] Furthermore, these maps are additive in the sense that given a second $K$-supported module $W_K$, we have \[\Add_{W_K} \circ \Add_{V_K} = \Add_{V_K \oplus W_K} \quad \text{and} \quad \MAdd_{W_K} \circ \MAdd_{V_K} = \MAdd_{V_K \oplus W_K} \; .\]
\end{proposition}
\begin{proof}
    First, we fix a point $\overline{V_K}$ in $\Rep(A,\dim(V_K))$ representing $V_K$, which amounts to choosing a basis for each $e_i V_K$. Consider the injective group homomorphism $\GL(v|_I) \hookrightarrow \GL(\hat{v}|_I)$  which comes from embedding \[\GL(v|_I) = \begin{pmatrix}\GL(v|_I) & 0 \\ 0 & 1_{\hat{v}_i - v_i} \end{pmatrix} \hookrightarrow \GL(\hat{v}|_I) \] for $i \in I$. With respect to this homomorphism, taking a direct sum with $V_K$ can be represented by equivariant closed embeddings of representation schemes \[\Add_{\overline{V_K}} \colon \Rep(Q,v) \hookrightarrow \Rep(Q,\hat{v}) \; , \quad \Rep(A,v) \hookrightarrow \Rep(A,\hat{v}) \; ,\] which can be defined as follows. The linear function on $\Rep(Q,\hat{v})$ which represents some $(i,j)$-matrix entry at some edge $h$ is
    \begin{itemize}
        \item sent to the corresponding linear function on $\Rep(Q,v)$ if $i \leq v_{t(h)}$ and $j \leq v_{s(h)}$,
        \item sent to $0$ if $i > v_{t(h)}$ and $j \leq v_{s(h)}$ or if $i \leq v_{t(h)}$ and $j > v_{s(h)}$, and
        \item sent to the value of the $(i-v_{t(h)},j-v_{s(h)})$-matrix entry, with respect to the chosen basis, of the same edge in $V_K$ if $i > v_{t(h)}$ and $j > v_{s(h)}$. Note that this is a constant.
    \end{itemize}
    We obtain the following commutative diagram of $\GL(v|_J)$-equivariant algebra homomorphisms:
    \[\begin{tikzcd}[row sep=small]
        \C[\Rep(A,\hat{v})] \arrow[r] & \C[\Rep(A,v)] \\
        \C[\Rep(A,\hat{v})]^{\GL(\hat{v}|_K)} \arrow[r] \arrow[u, phantom, sloped, "\subset"] & \C[\Rep(A,v)]^{\GL(v|_K)} \arrow[u, phantom, sloped, "\subset"]
    \end{tikzcd}\ .\]

    It is clear from the construction that the upper homomorphism is surjective. We will show that the same is true for the lower homomorphism by showing that the generators of $\C[\Rep(A,v)]^{\GL(v|_K)}$, as described by Theorem \ref{Thm:LBP-frozen}, lie in its image:
    \begin{itemize}
        \item Any generator of type (\ref{Thm:LBP-frozen:traces}), associated to some cycle in $Q$, is equal to the image of the generator in $\C[\Rep(A,\hat{v})]^{\GL(\hat{v}|_K)}$ corresponding to the same cycle, plus the trace of that cycle in $V_K$.
        \item Any generator of type (\ref{Thm:LBP-frozen:entries}), associated to some path $p$ in $Q$ and matrix coordinate $(i,j)$, where $1 \leq i \leq v_{t(p)}$ and $1 \leq j \leq v_{s(p)}$, equals the image of the generator of $\C[\Rep(A,\hat{v})]^{\GL(\hat{v}|_K)}$ corresponding to the same path and matrix coordinate.
    \end{itemize}
    Hence, we see that the bottom homomorphism is surjective as well, and taking its $\Spec$ we obtain a closed embedding of affine quotients $\Add_{V_K}$. We also see that this morphism is independent of the choice of representative $\overline{V_K}$, because the images of the generators of $\C[\Rep(A,\hat{v})]^{\GL(\hat{v}|_K)}$ are independent of the choice of basis.
    
    Clearly, $\Add_{V_K}$ is also $\GL(v|_J)$-equivariant. Now, by Lemma \ref{lem:successive-GIT}, the quiver schemes $\Mo_{I,\zeta}(A,v)$ and $\Mo_{I,\zeta}(A,\hat{v})$ are obtained by taking $\GL(v|_J)$-GIT quotients on both sides. Since reductive GIT quotients map closed invariant subschemes to closed subschemes, we also obtain the closed embedding of quiver schemes $\MAdd_{V_K}$.

    Finally, we show additivity. If we first pick representatives $\overline{V_K}$ and $\overline{W_K}$ by choosing bases of $V_K$ and $W_K$, then the union of the two bases will be a basis for the direct sum $V_K \oplus W_K$, which gives a representative $\overline{V_K \oplus W_K}$. With these choices, it is straight-forward to check, using the definition of $\Add$, that \[\Add_{\overline{W_K}} \circ \Add_{\overline{V_K}} = \Add_{\overline{V_K \oplus W_K}} \; .\] This equality is equivalent to an equality of the corresponding equivariant ring homomorphisms, which then also holds for maps between subspaces of invariants and relative invariants. This implies \[\Add_{W_K} \circ \Add_{V_K} = \Add_{V_K \oplus W_K} \quad \text{and} \quad \MAdd_{W_K} \circ \MAdd_{V_K} = \MAdd_{V_K \oplus W_K} \; .\]
\end{proof}

In the notation of Proposition \ref{prop:add-trivial}, note that any extension of $V$ by $V_K$, or of $V_K$ by $V$, will be represented by the same point as $V \oplus V_K$ in $\Rep(A,\hat{v}) \sslash \GL(\hat{v}|_K)$. Indeed, one can check that the invariant functions given by Theorem \ref{Thm:LBP-frozen} take the same value on any such extension as they do on the direct sum. Similarly, one can see that, given an extension \[0 \to V'_K \to V_K \to V''_K \to 0 \] of $A/(e_H)$-modules, we have equalities \[\Add_{V_K} = \Add_{V'_K \oplus V''_K} \quad \text{and} \quad \MAdd_{V_K} = \MAdd_{V'_K \oplus V''_K} \; .\]

\section{The Cornered Algebra}\label{sec:cornered}

\subsection{Morphisms induced by cornering}

\begin{assumptions}\label{asn:2}
    We fix all the assumptions and notation from Assumptions \ref{asn:1}. In addition, we assume that $J$ is non-empty, and that $\ida < \C Q$ is homogeneous with respect to the grading given by path length. Furthermore, we assume that the quotient algebra $A/(e_H)$ is isomorphic, as a graded algebra over $\C^K$, to the preprojective algebra $\Pi$ associated to a finite, possibly disconnected, ADE Dynkin diagram with vertex set $K$.
\end{assumptions}

\begin{remark}\label{rem:asn:2}
    The conditions of Assumptions \ref{asn:2} may look somewhat technical, but they are satisfied in some common settings, for instance when both the following conditions are fulfilled:
    
    \begin{enumerate}
        \item \label{rem:asn:2:1} $A$ is a preprojective algebra constructed from a finite undirected graph $G$, and 
        
        \item \label{rem:asn:2:2} removing the vertices of $H$ and all their adjacent edges from $G$ produces a (possibly disconnected) finite Dynkin diagram.
    \end{enumerate}
  
    In particular, Assumptions \ref{asn:2} are fulfilled when $A$ is the preprojective algebra $\Tilde{\Pi}$ for the \emph{affine} Dynkin diagram with vertex set $I$, with $J \subset I$ any nonempty subset. Other examples are given by the algebra $\Tilde{\Pi}^*$ of section \ref{sec:proof-rep-not-reduced} with $J=\{0\}$, $F=\emptyset$, or by the framed preprojective algebra $\hat{\Pi}$ and its quotient $A$ appearing in Remark \ref{rmk:CGGSisSuperfluous} and section \ref{sec:quot-nqv} with $F = \{\infty\}$, $J \subset I$ nonempty.
\end{remark}

By the assumption that $\ida$ is homogenous, $A$ inherits the grading by path length from $\C Q$. The length of a homogeneous element $a$ of $A$ will also be denoted by $l(a)$. 

We now define the \textit{cornered algebra} associated to $A$ and $H \subseteq Q_0$ as \[A_H = e_H A e_H \; .\] While $A_H \subset A$, note that $A_H$ is neither a quotient nor a subalgebra of $A$ in the unital sense: its unit is given by $e_H$. In order to define representation schemes for $A_H$, we need the following result.

\begin{lemma}\label{lem:finite-generation-cornered-algebra}
    $A_H$ is a finitely generated algebra.
\end{lemma}
\begin{proof}
    By Assumptions \ref{asn:2} and Theorem \ref{thm:finite-dimensional-finite-type}, $A/(e_H) \cong \Pi$ is a finite-dimensional algebra, graded by path length. Let $\lambda$ be the maximum length of a non-zero homogeneous element of $A/(e_H)$. Then any homogeneous element of $A$ of length $\geq \lambda + 1$ can be represented by a linear combination of paths in $\C Q$, each of which passes through $H$.

    Let $L \subset A_H$ be the image of the set of paths in $Q$ of length $\leq \lambda + 2$, starting and ending in $H$. This is a finite set; we claim that it generates $A_H$.

    It is sufficient to show that the class in $A_H$ of any path in $Q$ which starts and ends in $H$, can be generated by $L$. We proceed by induction on path length. To start, the image of any path of length $\leq \lambda+2$ is contained in the span of $L$ and therefore generated by it.
    
    For the induction step, suppose that for some $l \geq \lambda+2$ the class of every path of length $\leq l$, starting and ending in $H$, is generated by $L$. Let $p$ be a path of length $l+1$, starting and ending in $H$. Write $p=h_2p'h_1$ in $\C Q$, where $h_1$ and $h_2$ are edges with $s(h_1),t(h_2) \in H$. Now, if $t(h_1)=s(p')$ or $t(p')=s(h_2)$ is also in $H$, we are done because $p'$ and the respective edge will be generated from $L$ by the induction hypothesis. Hence, assume that $s(p'),t(p') \in K$. Then, as $l(p') = l-1 \geq \lambda+1$, $p'$ is equivalent modulo $\ida$ to a linear combination of paths that pass through $H$ at some point. But this implies that the image of $p$ in $A_H$ is equal to that of a linear combination of paths, each of which factors as a product of paths of shorter length which also start and end in $H$. The claim then follows from the induction hypothesis.
\end{proof}

We reiterate that Theorem \ref{thm:finite-dimensional-finite-type} (showing that $A/(e_H)$ is finite-dimensional) really is crucial for Lemma \ref{lem:finite-generation-cornered-algebra}. It is easy to construct examples for which $A/(e_H)$ is not finite-dimensional, and $A_H$ is not finitely generated: consider, for instance, a two-vertex quiver $Q$ given by

    \begin{center}\begin{tikzpicture}[ivert/.style={circle,draw=black,thick,inner sep=0pt,minimum size=5mm}]
    \node (i) at (0,0) [ivert] {$i$};
    \node (j) at (1.4,0) [ivert] {$j$};

    \draw[->] (i) to ["$x$", bend left=20] (j) ;
    \draw[->] (j) to ["$y$", bend left=20] (i);
    \draw[->] (j) to ["$z$"', loop] (j);
    \end{tikzpicture},\end{center}
and set $H = \{i\}$: we then have $\C Q/(e_H) \cong \C[z]$, and $e_H \C Q e_H\cong \C\langle\{yz^{k}x\}_{k>0}\rangle$, which is clearly not finitely generated.

We deduce from Lemma \ref{lem:finite-generation-cornered-algebra} that the cornered algebra can be written as the algebra of a quiver with relations:

\begin{corollary}\label{cor:finGenBetter}
     There exist a quiver $Q_H$ with vertex set $(Q_H)_0 = H$ and an ideal $\ida_H < \C Q_H$ such that \[A_H \cong \C Q_H / \ida_H \; .\]
\end{corollary}
\begin{proof}
    Consider a finite collection $L$ of generators of $A_H$, and assume that they are all homogeneous with respect to length, and source and target vertex. They span a vector space $\C L$, which is a finitely generated $\C^H$-bimodule generating $A_H$. Hence, there is a surjection of $\C^H$--algebras \[T^{\bullet}_{\C^H} (\C L) \twoheadrightarrow A_H \; ,\] where on the left-hand side we have take a tensor algebra of $\C L$ over $\C^H$. But this tensor algebra is nothing but the path algebra of a quiver with vertex set $H$ and set of edges $L$.
\end{proof}

\begin{remark}\label{rmk:CGGSisSuperfluous}
Consider an affine Dynkin diagram of type ADE, and adjoin a vertex $\infty$ and an edge connecting $\infty$ to the $0$-vertex. Let $\hat{Q}$ be the doubled quiver of this graph, and let $\hat{\Pi}$ be the preprojective algebra on $\hat{Q}$.

In \cite[Section 3]{craw2019punctual}, the authors introduce an algebra $A = \hat\Pi/(b^*)$, where $b^*$ is the unique arrow $0 \to \infty$. This is done in order to show \cite[Proposition 3.3]{craw2019punctual} that the cornered algebra $A_H$ is finitely generated (here, $H$ is a set of vertices that includes $\infty$ and at least one of the other vertices), which is used to construct a moduli space of $A_H$-modules. In any case, $A/(e_H)$ is certainly isomorphic to a preprojective algebra on a (possibly disconnected) finite Dynkin diagram, so Corollary \ref{cor:finGenBetter} is thus a generalisation of \cite[Proposition 3.3]{craw2019punctual}. Even more, it shows that one can in fact work with  $\hat\Pi_H$ directly, and there is no need to delete the arrow $b^*$.
\end{remark}

As a consequence of Corollary \ref{cor:finGenBetter}, we can construct moduli spaces of modules $\Mo_{J,\zeta_J}(A_H,v_H)$ by the methods of Section \ref{sec:quiver-varieties}. We have the following result.

\begin{theorem}\label{thm:cornered-embedding}
    Under Assumptions \ref{asn:2}, the following holds.
    \begin{enumerate}
        \item \label{thm:cornered-embedding:1} For $v \in \N^{Q_0}$ and $v_H=v|_H \in \N^H$, mapping $V \mapsto e_H V$ gives a closed, $\GL(v|_J)$-equivariant embedding of GIT quotients \[\Rep(A,v) \sslash \GL(v|_K) \hookrightarrow \Rep(A_H,v_H)\] and a closed embedding of quiver schemes \[\Mo_{I,\zeta}(A,v) \hookrightarrow \Mo_{J,\zeta_J}(A_H,v_H) \; .\]
        
        \item \label{thm:cornered-embedding:2} For $v \leq v' \in \N^{Q_0}$ with $v_H=v|_H=v'|_H$ the diagram \[\begin{tikzcd}\Mo_{I,\zeta}(A,v) \arrow[rd,hook] \arrow[r,hook] & \Mo_{I,\zeta}(A,v') \arrow[d,hook] \\ & \Mo_{J,\zeta_J}(A_H,v_H)\end{tikzcd}\] commutes, where the diagonal and vertical embedding are those of (\ref{thm:cornered-embedding:1}), and the horizontal embedding is given by adding modules supported on $K$, as constructed in Proposition \ref{prop:add-trivial}. 
        
        \item \label{thm:cornered-embedding:3}  For a fixed $v_H \in \N^H$, there is $v \in \N^{Q_0}$ such that $v|_H = v_H$, and for any $v' \geq v$ with $v'|_H = v_H$, the embedding of (\ref{thm:cornered-embedding:1}) is a bijection on closed points.
    \end{enumerate}
\end{theorem}

We give the proof in the next section. Let us also note that if the scheme $\Mo_{J,\zeta_J}(A_H,v_H)$ is reduced, the closed embedding inducing a bijection on closed points in Theorem \ref{thm:cornered-embedding}(\ref{thm:cornered-embedding:3}) will indeed be an isomorphism.

\begin{remark}
    It is worth remarking on the parallels between our Theorem \ref{thm:cornered-embedding} and Theorems 4.7, 1.4 of Craw--Yamagishi \cite{craw2026hilbertschemespointscanonical}. Again, in their case the algebra $A$ is defined as in Remark \ref{rmk:CGGSisSuperfluous}. In our language, they use $F=\{\infty\}$, $J=\{0\}$, and $K$ the vertex set of the finite preprojective algebra.
    
    At first, they follow the same strategy of taking successive GIT quotients, first with respect to $\GL(v|_K)$ (in their notation, $H_K$) and stability $0$, and then with respect to $\GL(v|_J)$ and general stability, to produce a closed immersion as in our Theorem \ref{thm:cornered-embedding}(\ref{thm:cornered-embedding:1}).
    
    However, at the expense of working with a restricted class of algebras and dimension vectors compared to our Assumptions \ref{asn:1}, \ref{asn:2}, they manage to obtain a bijection on closed points for a specific dimension vector of the form $v = (1,n\delta)$. This they deduce from the fact that, in the following diagram, the cornering morphism for stability $0$ is an isomorphism, the vertical VGIT morphisms are birational, and the target $\Mo_{J,\zeta_J}(A_H,(1,n))$ is known to be irreducible.
    \[\begin{tikzcd}
        \Mo_{I,\zeta}(A,(1,n\delta)) \arrow[r] \arrow[d] & \Mo_{J,\zeta_J}(A_H,(1,n))) \arrow[d]\\
        \Mo_{I,0}(A,(1,n\delta)) \arrow[r, "\sim"] & \Mo_{J,0}(A_H,(1,n))
    \end{tikzcd}\]
    They also identify the VGIT morphism on the right-hand side with the Hilbert--Chow-morphism \[\Hilb^n(\Sing) \to \Sym^n(\Sing) \; ,\] where $\Sing$ is the Kleinian singularity. We will use a generalisation of this identification in section \ref{sec:quot-nqv}.
\end{remark}

\subsection{Proof of Theorem \ref{thm:cornered-embedding}}\label{sec:proofs}

We start with the following technical lemma, whose proof is similar to that of Lemma \ref{lem:finite-generation-cornered-algebra}.

\begin{lemma}\label{lem:finite-generation-cornered-module}
    The $(A,A_H)$-bimodule $A e_H$ is finitely generated as a right $A_H$-module. The $(A_H,A)$-bimodule $e_H A$ is finitely generated as a left $A_H$-module.
\end{lemma}
\begin{proof}
        We show that $A e_H$ is a finitely generated right $A_H$-module. For the second part of the result, this proof can be read in a mirror. Note again that since $A/(e_H)$ is finite-dimensional, there is a maximal length for its elements, and denote this length by $\lambda$.

        Recall that $Q_0 = H\sqcup K$, which implies that we can write $Ae_H = (e_K+e_H)Ae_H = e_KAe_H \oplus A_H$. It will therefore suffice to show finite generation of $e_K A e_H$. Let $L \subset e_K A e_H$ be the set of classes of paths of length $\leq \lambda+1$ in $Q$ that start in $H$ and end in $K$. We claim that $L$ generates $e_K A e_H$. It suffices to show that the class in $e_K A e_H$ of any path which starts in $H$ and ends in $K$ can be generated over $A_H$ by $L$. We proceed by induction over the length of such paths. To start, the class of any path of length $\leq \lambda + 1$ is contained in the span of $L$ and therefore generated by it.

        For the induction step, suppose that for some $l \geq \lambda + 1$ the class of every path of length $\leq l$, starting in $H$ and ending in $K$, is generated over $A_H$ by $L$. Let $p \in e_K \C Q e_H$ be a path of length $l+1$. Write $p = p' h$ where $h$ is an edge and $p'$ a path of length $l$. If $s(p') = t(h) \in H$, then by the induction hypothesis the class of $p'$ in $e_K A e_H$ is generated over $A_H$ by $L$, and hence, so is $p$, because the class of $h$ is in $A_H$. Hence, we assume that $s(p') = t(h) \in K$. Then $l(p') = l \geq \lambda + 1$, so the image of $p'$ in $A/(e_H)$ is zero, so \[p' \equiv \sum_{i=1}^k \alpha_i p_i q_i \mod{\ida} \; , \] where $\alpha_i \in \C$, $p_i \in e_K \C Q e_H$ and $q_i \in e_H \C Q e_K$ for all $i$, and so \[p \equiv \sum_{i=1}^k \alpha_i p_i q_i h \mod{\ida} \; . \] For each $i$, the image of $q_i h$ is in $A_H$, the image of $p_i$ is in $e_K A e_H$, and $l(p_i) < l(p)$. The claim thus follows from the induction hypothesis.
\end{proof}

\begin{proof}[Proof of Theorem \ref{thm:cornered-embedding}]
    For part (\ref{thm:cornered-embedding:1}), we start with the observation that the map $V \mapsto e_H V$ induces a natural morphism of representation schemes \[\Rep(A,v) \rightarrow \Rep(A_H,v_H) \; .\] Indeed, this morphism can be constructed at the level of functors represented by the representation schemes (see Remark \ref{rem:rep-functor}), or using coordinates like in \cite[Prop.\ 4.4]{craw2026hilbertschemespointscanonical}. Furthermore, one can see in either description that the morphism is $\GL(v|_K)$--invariant and $\GL(v|_J)$--equivariant, hence corresponds to a $\GL(v|_J)$--equivariant algebra homomorphism \[\C[\Rep(A_H,v_H)] \rightarrow \C[\Rep(A,v)]^{\GL(v|_K)} \subset \C[\Rep(A,v)]\] (cf.\ \cite[Proof of Lemma 4.6]{craw2026hilbertschemespointscanonical} for equivariance). The image of this homomorphism is precisely the subalgebra generated by elements of type (\ref{Thm:LBP-frozen:entries}) in Theorem \ref{Thm:LBP-frozen} (where we have to replace the $I$ in the notation of Theorem \ref{Thm:LBP-frozen} by $K$). We will show that the elements of type (\ref{Thm:LBP-frozen:traces}) also lie in this image, and hence that the homomorphism is surjective.

    For this, notice that the element of type (\ref{Thm:LBP-frozen:traces}) associated to the class of a cycle is an invariant of the corresponding image of that cycle class in $A/[A,A]$ (notice that $[A,A]$ contains all the classes of paths whose source and target vertices do not coincide). Theorem \ref{thm:cocenter-finite-type} now gives \[\frac{A}{[A,A] + (e_H)} \cong \frac{\Pi}{[\Pi,\Pi]} \cong \C^K \; ,\] so the trace of any cycle of length $\geq 1$ is equal to a linear combination of traces of cycles that pass through $H$. These are also generated by elements of type (\ref{Thm:LBP-frozen:entries}).

    This concludes the proof of surjectivity of $\C[\Rep(A_H,v_H)] \rightarrow \C[\Rep(A,v)]^{\GL(v|_K)}$. Thus, we obtain a $\GL(v|_J)$--equivariant  closed embedding $\Rep(A,v) \sslash \GL(v|_K) \hookrightarrow \Rep(A_H,v_H)$. Taking the projective GIT quotient by $\GL(v|_J)$ now gives a closed embedding \[\Mo_{I,\zeta}(A,v) \hookrightarrow \Mo_{J,\zeta_J}(A_H,v_H) \; .\]

    Part (\ref{thm:cornered-embedding:2}) can easily be checked to hold at the level of representation schemes. The statement for GIT quotients follows immediately. 
    
    For part (\ref{thm:cornered-embedding:3}), consider an $A_H$-module $V_H$ with $\dim(V_H) = v_H$. Set \[V \coloneqq Ae_H \otimes_{A_H} V_H \; .\] Then $V$ is finite-dimensional by Lemma \ref{lem:finite-generation-cornered-module}, and $e_H V \cong V_H$. The same is true if we add to $V$ any finite number of vertex simple modules supported on $K$, increasing the $K$-components of the dimension vector of $V$. Note furthermore that we have the bound \[\dim_{\C}(V) \leq r \dim_{\C}({V_H}) \; ,\] where $r$ is the number of generators of $Ae_H$ as a right $A_H$-module. In particular, for fixed $v_H=\dim(V_H)$ the number of possible vectors $\dim(Ae_H \otimes_{A_H} V_H)$ is finite. If we pick $v$ to be their component-wise maximum, the above argument shows that we can write any $V_H$ of dimension $v_H$ as $V_H = e_HV$ for an $A$-module $V$ of dimension $v$. It follows that the map $\Rep(A,v) \rightarrow \Rep(A_H,v_H)$ induces a surjection on closed points, and so does the closed embedding $\Rep(A,v) \sslash \GL(v|_K) \hookrightarrow \Rep(A_H,v_H)$. But then this closed embedding is an isomorphism of underlying reduced subschemes, and so is its quotient $\Mo_{I,\zeta}(A,v) \hookrightarrow \Mo_{J,\zeta_J}(A_H,v_H)$.
\end{proof}

\begin{remark}\label{rem:is-bij-isomorphism}
    The inverse used in the proof of part (\ref{thm:cornered-embedding:3}) was constructed using the functor \[V_H \mapsto Ae_H \otimes_{A_H} V_H \; ,\] which is denoted by $j_!$ in \cite{craw2018multigraded}, and which is left adjoint to the functor $j^* \colon V \mapsto e_H V$. The functor $j_!$ does not, in general, preserve flat (that is, equidimensional) families of modules. Therefore, it seems unclear whether or how one could use this to define an actual inverse morphism to this bijection on closed points, which would establish a stronger version of Theorem \ref{thm:cornered-embedding}, (\ref{thm:cornered-embedding:3}).
    
    One approach to resolve this issue might be to consider the derived functor of $j_!$. Derived functors are used, for example, in \cite{Karmazyn17,SekiyaYamaura13} to construct isomorphisms between moduli schemes under certain conditions, but we do not know whether these ideas apply to the present setting.
\end{remark}

\begin{remark}
    We follow Remark \ref{rmk:CGGSisSuperfluous} by noting that it is also worthwhile to compare Theorem \ref{thm:cornered-embedding}(3) with \cite[Theorem 6.9]{CGGS2} and \cite[Theorem 4.8]{craw2019punctual}. Again, in the cited papers, the algebra $A$ is defined as in Remark \ref{rmk:CGGSisSuperfluous}, and $H$ is a vertex set consisting of $\infty$ and at least one other vertex.
    
    Let $v_H$ be a dimension vector for $A_H$-modules. In the proof of Theorem \ref{thm:cornered-embedding}, we give an upper bound $v$ for $\dim(Ae_H \otimes_{A_H} V_H)$, where $V_H$ ranges over all $A_H$-modules of dimension $v_H$. In the referenced theorem of \cite{craw2019punctual}, the authors argue for a slightly different result, namely that there is a \emph{lower} bound $v^\dagger$ for the set of dimension vectors $u$ such that every $A_H$-module $V_H$ of dimension vector $v_H$ can be written as $e_HV$ for some $A$-module $V$ of dimension $u$, in the specific case $v_H = n\delta|_H$. Indeed, our $v$ is simply the coordinate-wise maximum over over all vectors $\dim(Ae_H \otimes_{A_H} V_H)$, and one can construct examples for which $v>v^\dagger$.

    On the other hand, \cite[Theorem 6.9]{CGGS2}  is more directly comparable to Theorem \ref{thm:cornered-embedding}: in the proof of \emph{loc.cit} it is shown that, for \emph{any} dimension vector $v_H$ for $A_H$-modules, there is a dimension vector $v$ with $v|_H = v_H$, such that $\Mo_{I,\zeta}(A,v) \hookrightarrow \Mo_{J,\zeta_J}(A_H,v_H)$ is a bijection of the underlying reduced subschemes, just as we prove in Theorem \ref{thm:cornered-embedding}. The arguments used in \cite{CGGS2} to prove this depend on the fact that the algebras under consideration in \textit{ibid.} are preprojective algebras with one vertex deleted, and so Theorem \ref{thm:cornered-embedding} is a straightforward generalisation of this result to a larger class of algebras, and to other stability parameters. 
\end{remark}

\begin{remark}
    A careful inspection of the arguments presented in section \ref{sec:cornered} shows that Theorem \ref{thm:cornered-embedding} remains largely valid when weakening Assumptions \ref{asn:2} as follows. First, the only properties of $A/(e_H)$ that we use are
    \begin{enumerate}
        \item\label{asn:actual-1} $A/(e_H)$ is finite-dimensional, and
        \item\label{asn:actual-2} $A/(e_H) = [A/(e_H),A/(e_H)] \oplus \C^K$.
    \end{enumerate}
    Theorem \ref{thm:cornered-embedding} still holds if, instead of $A/(e_H)$ being a preprojective algebra, we only assume these properties (\ref{asn:actual-1}) and (\ref{asn:actual-2}).
    
    Second, the finite-dimensionality (\ref{asn:actual-1}) of $A/(e_H)$ is used only in the proofs of Lemmas \ref{lem:finite-generation-cornered-algebra} and \ref{lem:finite-generation-cornered-module}. If we assume the finite-generation properties stated in these lemmas to hold, we can further drop the assumption (\ref{asn:actual-1}) and Theorem \ref{thm:cornered-embedding} remains true.
\end{remark}

\section{Equivariant Quot Schemes as Nakajima Quiver Schemes}\label{sec:quot-nqv}

In this section, we use our results to give an independent proof of Theorem \ref{main_thm_appendix} below, another proof of which has appeared in \cite{craw2024orbifoldquotschemesle}. This theorem is concerned with a class of moduli spaces called \emph{orbifold Quot schemes}, originally introduced in \cite{CGGS2} and denoted by $\Quot_J^{v_J}([\Sing])$. Here, $\Sing$ is the Kleinian singularity associated to the Dynkin diagram via the McKay correspondence, and $[\Sing]$ is the corresponding orbifold, or DM stack. The definition of the Quot schemes is stated not in terms of $[\Sing]$, but rather in terms of the affine preprojective algebra $\Tilde{\Pi}$, the connection here being the Morita equivalence between $\Tilde{\Pi}$ and the skew-group algebra of $\Gamma$ acting on $\Aff^2$ \cite{reiten_van-den-bergh_two-dimensional-tame_89, crawley1998noncommutative}.

As in section \ref{sec:proof-rep-not-reduced}, let $I$ denote the vertex set of the affine Dynkin diagram, and let $0 \in I$ be the extended vertex. Here, we let $J \subseteq I$ be any subset.

\begin{definition}{\cite[Section 3.2]{CGGS2}}\label{def:OrbifoldQuotScheme}
    Choose a vector $v_J \in \mathbb N^J$. The scheme $\Quot_J^{v_J}([\Sing])$ is the fine moduli space of $\tilde\Pi_J$-module surjections $ e_J\tilde\Pi e_0\to Z$ such that $\dim(Z) = v_J$. 
\end{definition}

On the other hand, we shall denote by $\NQV_{\zeta}(v,w)$ the Nakajima quiver scheme associated to the affine Dynkin diagram, stability parameter $\zeta$, dimension vector $v$, and framing vector $w$. In this section, the only framing vector appearing will be $w = \Lambda_0 \coloneqq (1,0,\ldots,0)$, i.e., a one-dimensional framing attached to the vertex $0$. As in Remark \ref{rmk:CGGSisSuperfluous}, we let $\hat{\Pi}$ denote the preprojective algebra associated to the Dynkin graph framed with a vertex $\infty$ at $0$. We let $b$ be the class of the edge in $\hat{\Pi}$ going from $\infty$ to $0$ and let $b^*$ be its opposite. We also define \[A \coloneqq \hat{\Pi}/(b^*) \; .\] The quivers generating the algebras $\Tilde{\Pi}$, $A$, and $\hat{\Pi}$ are depicted below.

\vspace{6mm}

\begin{center}
\begin{tikzpicture}[framing/.style={rectangle,draw=black,thick,inner sep=0pt,minimum size=5mm},ivert/.style={circle,draw=black,thick,inner sep=0pt,minimum size=5mm}]
    \node (a) at (0,0) {$\Tilde{\Pi}$};
    \node (tl) at (3,.8) [ivert] {$0$};
    \node (bl) at (3,-.8) [ivert] {};
    \node (l) at (4,0) [ivert] {};
    \node (r) at (5.3,0) [ivert] {};
    \node (tr) at (6.3,.8) [ivert] {};
    \node (br) at (6.3,-.8) [ivert] {};

    \draw[->] (tl) to [bend left=20] (l);
    \draw[->] (l) to [bend left=20] (tl);
    \draw[->] (bl) to [bend left=20] (l);
    \draw[->] (l) to [bend left=20] (bl);
    \draw[->] (l) to [bend left=20] (r);
    \draw[->] (r) to [bend left=20] (l);
    \draw[->] (tr) to [bend left=20] (r);
    \draw[->] (r) to [bend left=20] (tr);
    \draw[->] (br) to [bend left=20] (r);
    \draw[->] (r) to [bend left=20] (br);
\end{tikzpicture}
\\ \vspace{6mm}
\begin{tikzpicture}[framing/.style={rectangle,draw=black,thick,inner sep=0pt,minimum size=5mm},ivert/.style={circle,draw=black,thick,inner sep=0pt,minimum size=5mm}]
    \node (a) at (0,0) {$A$};
    \node (f) at (1.7,.8) [framing] {$\infty$};
    \node (tl) at (3,.8) [ivert] {$0$};
    \node (bl) at (3,-.8) [ivert] {};
    \node (l) at (4,0) [ivert] {};
    \node (r) at (5.3,0) [ivert] {};
    \node (tr) at (6.3,.8) [ivert] {};
    \node (br) at (6.3,-.8) [ivert] {};

    \draw[->] (f) to ["$b$"] (tl);
    \draw[->] (tl) to [bend left=20] (l);
    \draw[->] (l) to [bend left=20] (tl);
    \draw[->] (bl) to [bend left=20] (l);
    \draw[->] (l) to [bend left=20] (bl);
    \draw[->] (l) to [bend left=20] (r);
    \draw[->] (r) to [bend left=20] (l);
    \draw[->] (tr) to [bend left=20] (r);
    \draw[->] (r) to [bend left=20] (tr);
    \draw[->] (br) to [bend left=20] (r);
    \draw[->] (r) to [bend left=20] (br);
\end{tikzpicture}
\\ \vspace{6mm}
\begin{tikzpicture}[framing/.style={rectangle,draw=black,thick,inner sep=0pt,minimum size=5mm},ivert/.style={circle,draw=black,thick,inner sep=0pt,minimum size=5mm}]
    \node (a) at (0,0) {$\hat{\Pi}$};
    \node (f) at (1.7,.8) [framing] {$\infty$};
    \node (tl) at (3,.8) [ivert] {$0$};
    \node (bl) at (3,-.8) [ivert] {};
    \node (l) at (4,0) [ivert] {};
    \node (r) at (5.3,0) [ivert] {};
    \node (tr) at (6.3,.8) [ivert] {};
    \node (br) at (6.3,-.8) [ivert] {};

    \draw[->] (f) to ["$b$", bend left=20] (tl);
    \draw[->] (tl) to ["$b^*$", bend left=20] (f);
    \draw[->] (tl) to [bend left=20] (l);
    \draw[->] (l) to [bend left=20] (tl);
    \draw[->] (bl) to [bend left=20] (l);
    \draw[->] (l) to [bend left=20] (bl);
    \draw[->] (l) to [bend left=20] (r);
    \draw[->] (r) to [bend left=20] (l);
    \draw[->] (tr) to [bend left=20] (r);
    \draw[->] (r) to [bend left=20] (tr);
    \draw[->] (br) to [bend left=20] (r);
    \draw[->] (r) to [bend left=20] (br);
\end{tikzpicture}
\end{center}

\vspace{6mm}

We can now construct our Nakajima quiver scheme as \[\NQV_{\zeta}(v,\Lambda_0) = \Mo_{I, \zeta}(\hat{\Pi},(1,v)) \; .\] (Note that the definition of Nakajima quiver schemes generally does not include any relations at the framing vertices. In the case of one-dimension framing at $0$, however, the additional relation $b^*b=0$ follows from the preprojective relations at the other vertices. \cite[Lemma 25]{kuznetsov2007quiver})

\begin{theorem}\label{main_thm_appendix}
    Let $J \subseteq I$ be a non-empty subset and $v_J \in \N^J$ be a dimension vector. The orbifold Quot scheme $\Quot_J^{v_J} ([\Sing])$ is non-empty if and only if the Nakajima quiver scheme $\mathfrak{M}_{\zeta^{\bullet}}(v,\Lambda_0)$ is non-empty for some vector $v\in \N^I$ satisfying $v|_J=v_J$. In this case, for any sufficiently large choice of $v \in \N^I$ with $v|_J=v_J$, we have an isomorphism \[\Quot_J^{v_J} ([\Sing])_{\rm red} \cong  \mathfrak{M}_{\zeta^{\bullet}}(v,\Lambda_0)_{\rm red},\] where on both sides we take the reduced scheme structure.
\end{theorem}

Before we prove Theorem \ref{main_thm_appendix}, we recall the following result. While this fact is well known to experts, we were unable to find a complete proof in this generality, so we will give one below. (See also \cite[Lemma 3.1]{craw2019punctual} or \cite[Lemma 35]{kuznetsov2007quiver} for related arguments.)

\begin{lemma}\label{lem:j-zero}
    Let $\zeta \in \N^I$ be such that $\zeta>0$, and let $v \in \N^I$ be arbitrary. Let $V$ be a $\hat{\Pi}$-module of dimension vector $(1,v)$ that is $\zeta$-stable as a framed $\Tilde{\Pi}$-module. Then $b^*z=0$ for all $z\in V$.
\end{lemma}
\begin{proof}
    Let $\Tilde{Q}$ be the affine double quiver generating $\Tilde{\Pi}$. Crawley-Boevey \cite{crawley1998noncommutative} describes a Morita equivalence between \[\C \Tilde{Q} \quad \text{and} \quad \C\langle x,y\rangle \rtimes \Gamma \; .\] Here, $\C\langle x,y\rangle \rtimes \Gamma$ is the skew-group algebra associated to the action of $\Gamma$ on the free associative algebra $\C\langle x,y\rangle$, and the vertex idempotents on the left correspond to the idempotents in $\C \Gamma$ associated to the irreducible representations of $\Gamma$ on the right. Under this Morita equivalence, the preprojective relations $p=0$ on the left-hand side correspond to the relation $[x,y]=0$ on the right-hand side, which gives rise to the Morita equivalence between \[\Tilde{\Pi} \quad \text{and} \quad \C[x,y] \rtimes \Gamma \; .\]

    Clearly, a $\hat{\Pi}$-module of dimension vector $(1,v)$ is the same thing as a $\C \Tilde{Q}$-module $V$ of dimension vector $v$ together with linear maps $b \colon \C \to V_0$ and $b^* \colon V_0 \to \C$ such that \[p z = b b^* e_0 z \quad \text{for all} \; z \in V \; .\] Under the Morita equivalence above, this now corresponds to a $\C\langle x,y\rangle \rtimes \Gamma$-module $U$ with maps $a \colon \C \to U^{\Gamma}$ and $a^* \colon U^{\Gamma} \to \C$ so that \[[x,y]z' = a a^* e_0 z' \quad \text{for all} \; z' \in U \; ,\] where $e_0 = \sum_{\gamma \in \Gamma} \gamma$ is now the $\Gamma$-invariant idempotent.
    
    The condition that $V$ is $\zeta$-stable is equivalent to the condition that $V$ is generated as a module by the image of $b$. That is equivalent to the condition that the Morita partner $U$ is generated by the image of $a$.

    By construction of the skew-group algebra, a $\C\langle x,y\rangle \rtimes \Gamma$-module $U$ is the same thing as a $\Gamma$-representation $U$ together with endomorphisms $x,y \in \End(U)$ so that the $\Gamma$-action on $\C x \oplus \C y$ is compatible with that on $U$. Using this formulation, we can now forget about the $\Gamma$-action and apply \cite[Proposition 2.8]{nakajima1999lectures} to see that $a^*=0$. This implies that $b^*=0$ on the Morita partner $V$.
\end{proof}

\begin{proof}[Proof of Theorem \ref{main_thm_appendix}]
It was shown in \cite[Prop.\ 4.2]{CGGS2} that there is an isomorphism \[\Quot_J^{v_J} ([\Sing]) \cong \Mo_{\zeta_J}(A_H,(1,v_J)) \; .\] On the other hand, Theorem \ref{thm:cornered-embedding} shows that there is a closed embedding \[\Mo_{\zeta^{\bullet}}(A,(1,v)) \hookrightarrow \Mo_{\zeta_J}(A_H,(1,v_J))\] which, possibly after increasing $v$ on $K$, becomes a bijection on closed points, and hence an isomorphism of underlying reduced schemes \[\Mo_{\zeta^{\bullet}}(A,(1,v))_{\rm red} \cong \Mo_{\zeta_J}(A_H,(1,v_J))_{\rm red} \; .\]
    
Finally, for any $v \in \N^I$, we get a closed embedding \[\Mo_{\zeta^{\bullet}}(A,(1,v)) \hookrightarrow \NQV_{\zeta^{\bullet}}(v,\Lambda_0) \, \] coming from the fact that $A$ is a quotient of the framed preprojective algebra defining $\NQV_{\zeta^{\bullet}}(v,\Lambda_0)$, which has an additional arrow $0 \to \infty$. Recall from \cite{NakajimaBranching} that for this choice of framing, a $\zeta^\bullet$-stable module is either a one-dimensional module $S_i$ supported at a vertex $i\in K$, or contains the one-dimensional framing. Now pick a closed point in $\NQV_{\zeta^{\bullet}}(v,\Lambda_0)$ represented by a polystable module $V = V' \oplus \bigoplus_{i \in K} S_i^{\oplus v^s_i}$ where $V'$ is $\zeta^{\bullet}$-stable. Then $V'$ is also $\zeta$-stable for any generic stability condition $\zeta>0$. By Lemma \ref{lem:j-zero}, $0 \to \infty$ acts as zero on the stable module $V'$, 
and hence on $V$. The point represented by $V$ therefore lies in the image of $\Mo_{\zeta^{\bullet}}(A,(1,v)) \hookrightarrow \NQV_{\zeta^{\bullet}}(v,\Lambda_0)$. This shows that this embedding is surjective on closed points and hence we have an isomorphism of reduced subschemes \[\Mo_{\zeta^{\bullet}}(A,(1,v))_{\rm red} \cong \NQV_{\zeta^{\bullet}}(v,\Lambda_0)_{\rm red} \; .\] Combining these three isomorphisms, we deduce Theorem~\ref{main_thm_appendix}.
\end{proof}

The vector $v$ used in this proof is given by Theorem \ref{thm:cornered-embedding} (\ref{thm:cornered-embedding:3}). In \cite[Prop.\ 3.6]{BGScollapsing-fibers}, another argument for the existence of such $v$ is given, but for the corresponding Nakajima quiver schemes (the notation used there being $v^{\oplus}$).

\bibliographystyle{amsplain}
\bibliography{References}

\end{document}